\documentclass[11pt]{amsart}
\usepackage{amsfonts,amsthm,amsmath}
\usepackage{natbib}
\usepackage{hyperref}

\newtheorem{theorem}{Theorem}[section]
\newtheorem{lemma}{Lemma}[section]

\newcommand{\weak}{\stackrel{w}{\longrightarrow}}
\newcommand{\prob}{\stackrel{P}{\longrightarrow}}
\newcommand{\eid}{\stackrel{d}{=}}
\newcommand{\one}{{\bf 1}}
\newcommand{\reals}{{\mathbb R}}
\newcommand{\bbr}{\reals}
\newcommand{\vep}{\varepsilon}

\numberwithin{equation}{section}

\begin{document}
\title[Truncated heavy tails]{Central Limit Theorem for truncated heavy tailed Banach valued random vectors}
\author[A. Chakrabarty]{Arijit Chakrabarty}
\address{Department of Mathematics,
Indian Institute of Science,
Bengaluru - 560012, India}
\email{arijit@math.iisc.ernet.in}

\begin{abstract}
In this paper the question of the extent to which truncated heavy
tailed random vectors, taking values in a Banach space, retain the characteristic features of heavy
tailed random vectors, is answered from the point of view of the
central limit theorem.
\end{abstract}
 \keywords{ heavy
tails, truncation, regular variation, central limit
  theorem, probability on Banach spaces\vspace{.5ex}}
\thanks{Research partly supported by the NSF grant ``Graduate and Postdoctoral Training in Probability and its Applications'' at
Cornell University, and the Centenary Postdoctoral Fellowship at the Indian Institute of Science.}
\maketitle

\section{Introduction} \label{sec:intro}
Situations where heavy-tailed distributions is a good fit, and at the same time there is a physical upper bound on the quantity of interest, are common in nature. Clearly, the natural model for phenomena like this is a
truncated heavy-tailed distribution - a distribution that matches
a heavy-tailed one till a specified limit and after that it decays
significantly faster or simply vanishes.
This leads to the general
question: when can the upper bound be considered to be large
enough so that the effect of truncating by that is negligible? The first attempt at answering this question, in finite dimensional spaces, was made in \cite{chakrabarty:samorodnitsky:2009}. In the current paper, the investigation started by that paper has been continued to achieve similar results in Banach spaces.

Suppose that $B$ is a separable Banach space and that $H, H_1, H_2,\ldots$ are $B$-valued random variables in the domain of attraction of an $\alpha$-stable random variable ${\mathcal V}$ with $0<\alpha<2$. This means that there are sequences $a_n$ and $b_n$ so that as $n\longrightarrow\infty$,
\begin{equation}\label{DOA}
b_n^{-1}\left(\sum_{j=1}^nH_j-a_n\right)\Longrightarrow{\mathcal V}\,.
\end{equation}
We assume that the truncating threshold goes to infinity along with the sample size, and hence we essentially have a sequence of models. We denote both - the sample size and the number of the model by $n$, and the truncating threshold in the $n$-th model by $M_n$. The $n$th row of the
triangular array will consist of observations $X_{nj},\,
j=1,\ldots, n$, which are assumed to be generated according to the
following mechanism:
\begin{equation}\label{e:the.model}
X_{nj}:=H_j\one\left(\|H_j\|\le  M_n\right)
+\frac{H_j}{\|H_j\|}(M_n+L_j)\one\left(\|H_j\|>M_n\right)\,,
\end{equation}
$j=1,\ldots, n, \, n=1,2,\ldots$. Here $(L,L_1,L_2,\ldots)$ is a sequence of i.i.d. nonnegative random variables independent of $(H,H_1,H_2,\ldots)$. For each $n=1,2,\ldots$ we view the observation $X_{nj}$, $j=1,\ldots,n$ as having power tails that are truncated at level $M_n$. The random variable $L$ can be thought of as to model that outside the ball of radius $M_n$, the tail ``decays significantly faster or simply vanishes''. $L$ is assumed to have finite second moment.

In \cite{chakrabarty:samorodnitsky:2009} two regimes depending on the growth rate of $M_n$ and the tail of the random variable $H$ were introduced as follows: the tails in
the model \eqref{e:the.model} are said to be
\begin{equation} \label{e:regimes}
\begin{array}{ll}
\text{truncated softly} & \text{if} \ \lim_{n\to\infty}nP\left( \|
H\|>M_n\right) = 0\,,\\
\text{truncated hard} & \text{if} \ \lim_{n\to\infty}nP\left( \|
H\|>M_n\right) = \infty\,.
\end{array}
\end{equation}
It was shown in that paper that as far as the central limit behavior of the row sum is concerned, observations with
softly truncated tails behave like heavy tailed random variables,
while observations with hard truncated tails behave like light
tailed random variables. In Theorem \ref{clt.t2}, the main result of this paper, we show that the result under hard truncation can be extended to Banach spaces, if the ``small ball criterion'' holds. Doing this is not straightforward because of the following reason. While in finite-dimensional spaces, convergence in law is equivalent to one-dimensional convergence of each linear functional in law to the linear functional evaluated at the limit, the same is not true in Banach spaces. In the latter spaces, one needs to check in addition some tightness conditions; see for example, \cite{ledoux:talagrand:1991} or \cite{araujo:gine:1980} for details.

Section \ref{sec:clt} contains the results and their proofs. A couple of examples are studied in Section \ref{sec:examples} - one where the hypothesis of Theorem \ref{clt.t2} can be checked, and the other where the claim of that result does not hold. The examples serve the purpose of showing that there is a need for such a result, and that the result has some practical value.

\section{A Central Limit Theorem for truncated heavy-tailed random variables}\label{sec:clt}
The triangular array $\{X_{nj}:1\le j\le n\}$ is as defined in \eqref{e:the.model}. We would like to know if the row sums $S_n$, defined by
\begin{equation}\label{rowsum}
S_n:=\sum_{j=1}^nX_{nj}\,,
\end{equation}
still converge in law after appropriate centering and scaling. Exactly same arguments as those in the proof of Theorem 2.1 in \cite{chakrabarty:samorodnitsky:2009} show that if the truncated heavy-tailed model is in the soft truncation regime as defined in \eqref{e:regimes}, then
$$
b_n^{-1}(S_n-a_n)\Longrightarrow{\mathcal V}\,.
$$
In other words, from the point of view of central limit behavior
of the partial sums, the truncated heavy-tailed model retains much
of the heavy-tailedness. Hence, we shall assume
throughout that the model is in the hard truncation regime, {\it
i.e.},
\begin{equation}\label{ht}
\lim_{n\to\infty}nP(\|H\|>M_n)=\infty\,.
\end{equation}

As mentioned earlier, easy-to-check criteria for  satisfying the Central Limit Theorem on Banach spaces are not known. An example of
the not-so-easy-to-check  ones is Theorem
10.13, page 289 in \cite{ledoux:talagrand:1991}, known as the ``small ball criterion''.
The main result of this paper, Theorem \ref{clt.t2}, is an analogue of this theorem in the
truncated setting under hard truncation. But before stating that, we need the following preliminary.
It is known that \eqref{DOA} implies that there is a probability measure $\sigma$ on
$$
{\mathcal S}:=\{x\in B:\|x\|=1\}
$$
such that as $t\longrightarrow\infty$,
\begin{equation}\label{ld.spectral}
{P\left(\frac H{\|H\|}\in\cdot\bigg|\|H\|>t\right)}\weak\sigma(\cdot)
\end{equation}
weakly on $\mathcal S$; see Corollary 6.20(b), page 151 in \cite{araujo:gine:1980}.

\begin{theorem}\label{clt.t2} There is a Gaussian measure $\gamma$ on $B$ such that
\begin{equation}\label{clt.eq2}
B_n^{-1}(S_n-ES_n)\Rightarrow\gamma
\end{equation}
if and only if the following hold:
\begin{enumerate}
\item (small ball criterion) For every $\epsilon>0$
$$\liminf_{n\rightarrow\infty}P(B_n^{-1}\|S_n-ES_n\|<\epsilon)>0\,,$$
\item $\sup_{n\ge1}B_n^{-1}E\|S_n-ES_n\|<\infty$,
\end{enumerate}
where
$$
B_n:=\left[nM_n^2P(\|H\|>M_n)\right]^{1/2}\,.
$$
In that case, the characteristic function of $\gamma$ is given by
\begin{equation}\label{clt.e3}
\hat\gamma(f)=\exp\left(-\frac2{2-\alpha}\int_Sf^2(s)\sigma
(ds)\right)\,,f\in B^\prime\,.
\end{equation}
Here, $B^\prime$ is the dual of $B$, the space of linear functionals on $B$.
\end{theorem}

For the proof, we shall need the following one-dimensional lemma, which follows by exactly similar arguments as those in Theorem 2.2 of \cite{chakrabarty:samorodnitsky:2009}, and hence we omit the proof.

\begin{lemma}\label{clt.l1}For every $f$ in $B^\prime$,
$$B_n^{-1}(f(S_n)-Ef(S_n))\Rightarrow N\left(0,\frac2{2-\alpha}\int_Sf^2(s)\sigma (ds)\right)\,.$$
\end{lemma}

\begin{proof}[Proof of Theorem \ref{clt.t2}] First we prove the direct part, {\it i.e.}, we assume that 1. and 2. hold. We first show that it suffices to check that $\{{\mathcal L}(Z_n)\}$ is relatively compact where
\begin{eqnarray*}
Z_n&:=&B_n^{-1}\sum_{j=1}^nY_{nj}\,,\\
Y_{nj}&:=&X_{nj}-X_{nj}^\prime\,
\end{eqnarray*}
and for every $n$, $X_{n1}^\prime,X_{n2}^\prime,\ldots$ are i.i.d.
copies of $X_{n1}$ so that $(X_{nj}^\prime:j\ge1)$  and
$(X_{nj}:j\ge1)$ are independent families. To see this, suppose that we have shown that $\{{\mathcal L}(Z_n)\}$ is relatively compact.
By Corollary 4.11, page 27 in \cite{araujo:gine:1980}, it follows that the sequence
$\{{\mathcal L}(B_n^{-1}S_n)\}$ is relatively shift compact, {\it
i.e.}, there exists some sequence $\{v_n\}\subset B$ such that
$\{{\mathcal L}(B_n^{-1}S_n-v_n)\}$ is relatively compact. By Theorem
4.1 in de Acosta and Gin\'e (1979), for relative compactness of $\{{\mathcal L}[B_n^{-1}(S_n-ES_n)]\}$, it suffices to check that
\begin{equation}\label{clt.t2.eq7}
\lim_{t\to\infty}\limsup_{n\to\infty}nE\left[\|U_n\|\one(\|U_n\|>t)\right]=0\,,
\end{equation}
where
\begin{equation}\label{defU}
U_n:=B_n^{-1}\left[H\one(\|H\|\le M_n)+\frac{H}{\|H\|}(M_n+L)\one(\|H\|>M_n))\right]\,.
\end{equation}
By \eqref{ht}, it follows that $B_n\gg M_n$. Thus, for fixed $t>0$ and $n$ large enough,
$$
nE\left[\|U_n\|\one(\|U_n\|>t)\right]
$$
$$
=nB_n^{-1}P(\|H\|>M_n)\left\{M_nP(L>B_nt-M_n)+E\left[L\one(L>B_nt-M_n)\right]\right\}\,.
$$
Since $EL^2<\infty$,
$$
E\left[L\one(L>B_nt-M_n)\right]\le\frac{E(L^2)}{B_nt-M_n}
$$
and
$$
M_nP(L>B_nt-M_n)\le \frac{M_n}{(B_nt-M_n)^2}EL^2=o(B_n^{-1})
$$
as $n\longrightarrow\infty$. Thus, for all fixed $t>0$,
\begin{eqnarray}
\lim_{n\to\infty}nE\left[\|U_n\|\one(\|U_n\|>t)\right]
=0\,.\label{new1}
\end{eqnarray}
This shows \eqref{clt.t2.eq7} and hence that $\{{\mathcal L}[B_n^{-1}(S_n-ES_n)]\}$ is relatively compact. In view of Lemma \ref{clt.l1}, this will complete the proof of the direct part.

First we record some
properties of the random variables defined above, which shall be
used in the proof. The hypotheses immediately imply that for all
$\epsilon>0$
\begin{equation}\label{clt.t2.eq1}
\liminf_{n\rightarrow\infty}P(\|Z_n\|<\epsilon)>0
\end{equation}
and that
\begin{equation}\label{clt.t2.eq2}
\sup_{n\ge1}E\|Z_n\|<\infty\,.
\end{equation}
Let $\{F_k\}$ be any sequence of increasing finite-dimensional
subspaces so that
\begin{equation}\label{clt.closure}
\mbox{closure}\left(\bigcup_{k=1}^\infty F_k\right)=B\,.
\end{equation}
For any subspace $F$ of $B$, denote by $T_F$ the canonical map from $B$ to the quotient space $B/F$. By Corollary 6.19 (page 151) in \cite{araujo:gine:1980}, it
follows that for every $k$, $T_{F_k}(H)$ is in the domain of attraction of some
$\alpha$-stable law with the same scaling constant $(b_n)$ as that
of $H$, and that
\begin{equation}\label{clt.t2.eq3}
\lim_{k\rightarrow\infty}\sup_{n\ge1}nP\left(\|T_{F_k}(H)\|>b_n\right)=0\,.
\end{equation}
Clearly, for every $k$, there is $C_k\in[0,\infty)$ so that as
$t\longrightarrow\infty$, $$P(\|T_{F_k}(H)\|>t)\sim
C_kP(\|H\|>t)\,.$$ It follows by (\ref{clt.t2.eq3}) that
$\lim_{k\rightarrow\infty}C_k=0$.
Note that,
\begin{eqnarray*}
&&E\|T_{F_k}(X_{n1})\|^2\\
&=&E\left[\|T_{F_k}(H)\|^2\one(\|H\|\le M_n)\right]\\
&&\,\,\,\,+E\left[\frac{\|T_{F_k}(H)\|^2}{\|H\|^2}(M_n+L)^2\one(\|H\|>M_n)\right]\\
&\le&E\left[\|T_{F_k}(H)\|^2\one\left(\|T_{F_k}(H)\|\le M_n\right)\right]\\
&&\,\,\,\,+E\left[\frac{\|T_{F_k}(H)\|^2}{\|H\|^2}\one(\|H\|>M_n)\right]E(M_n+L)^2\,.
\end{eqnarray*}
By the Karamata theorem (Theorem B.1.5, page 363 in \cite{dehaan:ferreira:2006}),
$$\lim_{n\rightarrow\infty}[M_n^2P(\|H\|>M_n)]^{-1}E\left(\|T_{F_k}(H)\|^2\one\left(\|T_{F_k}(H)\|\le M_n\right)\right)=\frac\alpha{2-\alpha}C_k\,.$$
By \eqref{ld.spectral}, it follows that as $n\longrightarrow\infty$,
$$E\left[\frac{\|T_{F_k}(H)\|^2}{\|H\|^2}\one(\|H\|>M_n)\right]\sim P(\|H\|>M_n)\int_{\mathcal S}\|T_{F_k}(s)\|^2\sigma(ds)\,.$$
That \eqref{clt.closure} holds and the fact that $\sigma$ is a finite measure implies that
$$
\lim_{k\to\infty}\int_{\mathcal S}\|T_{F_k}(s)\|^2\sigma(ds)=0\,.
$$
Thus, in view of the assumption that $EL^2<\infty$, it follows that
\begin{equation}\label{clt.new}
\lim_{k\rightarrow\infty}\limsup_{n\rightarrow\infty}[M_n^2P(\|H\|>M_n)]^{-1}E\|T_{F_k}(X_{n1})\|^2=0\,,
\end{equation}
which in turn implies that
\begin{equation}\label{clt.t2.eq4}
\lim_{k\rightarrow\infty}\limsup_{n\rightarrow\infty}[M_n^2P(\|H\|>M_n)]^{-1}E\|T_{F_k}(Y_{n1})\|^2=0\,.
\end{equation}

Coming to the proof, in view of the criterion for relative
compactness discussed in \cite{ledoux:talagrand:1991} (page
40-41), it suffices to show that given $ \epsilon>0$, there is a
finite dimensional subspace $F$ with
\begin{equation}\label{clt.tight}
\limsup_{n\rightarrow\infty}P\left[\|T_F(Z_n)\|> \epsilon\right]\le \epsilon\,.
\end{equation}
Let $ \vep_1, \vep_2,\ldots$ be an i.i.d. sequence of
Rademacher random variables, independent of
$(X_n,X_n^\prime,n\ge1)$, and let $E_ \vep$ denote the
conditional expectation given $\{Y_{nj}\}$.  It suffices to show that for all $\eta>0$,
$$
\lim_{k\rightarrow\infty}\limsup_{n\rightarrow\infty}P\Biggl[\Biggl|\biggl\|\sum_{j=1}^n \vep_jT_{F_k}(Y_{nj})\biggr\|
$$
\begin{equation}\label{clt.t2.eq5}
-E_ \vep\biggl\|\sum_{j=1}^n \vep_jT_{F_k}(Y_{nj})\biggr\|\Biggr|>B_n\eta\Biggr]=0\,,
\end{equation}
 and that there is a numerical constant $C>0$ so that for every
$\delta>0$,
\begin{equation}\label{clt.t2.eq6}
\limsup_{k\rightarrow\infty}\limsup_{n\rightarrow\infty}P\left[E_ \vep\biggl\|\sum_{j=1}^n \vep_jT_{F_k}(Y_{nj})\biggr\|>B_nC\delta\right]<\delta\,,
\end{equation}
whenever $\{F_k\}$ is an increasing sequence of finite-dimensional subspaces satisfying \eqref{clt.closure}.

To establish (\ref{clt.t2.eq5}), it suffices to check that
$$
\lim_{k\rightarrow\infty}\limsup_{n\rightarrow\infty}P\Biggl[\Biggl|\biggl\|\sum_{j=1}^n \vep_jT_{F_k}(u_{nj})\biggr\|
$$
\begin{equation}\label{clt.t2.trunc}
\,\,\,\,\,\,\,-E_ \vep\biggl\|\sum_{j=1}^n \vep_jT_{F_k}(u_{nj})\biggr\|\Biggr|>B_n\eta\Biggr]=0
\end{equation}
where
$$u_{nj}:=Y_{nj}\one{\left(\|Y_{nj}\|\le \beta B_n\right)}\,,$$
$\beta >0$ is to be specified later. This is because for $n$ large
enough,
\begin{eqnarray*}
&&B_n^{-1}E\biggl\|\sum_{j=1}^nY_{nj}\one\left(\|Y_{nj}\|>\beta B_n\right)\biggr\|\\
&\le&nB_n^{-1}E\left[\|Y_{n1}\|\one(\|Y_{n1}\|>\beta B_n)\right]\\
&\le&nB_n^{-1}E\Biggl[\left(\|X_{n1}\|+\|X^\prime_{n1}\|\right)\biggl\{\one\left(\|X_{n1}\|>\frac\beta2 B_n\right)\\
&&\,\,\,\,\,\,\,\,\,\,\,\,\,\,\,\,\,\,\,\,+\one\left(\|X^\prime_{n1}\|>\frac\beta2 B_n\right)\biggr\}\Bigg]
\end{eqnarray*}
\begin{eqnarray*}
&=&2nB_n^{-1}E\left[(M_n+L)\one(\|H\|>M_n)\one\left(L>\frac\beta2B_n-M_n\right)\right]\\
&&\,\,\,\,\,+2nB_n^{-1}P(\|H\|>M_n)P\left(L>\frac\beta2B_n-M_n\right)E\|X_{n1}\|\\
&=:&Q_1+Q_2\,.
\end{eqnarray*}
Clearly,
\begin{eqnarray*}
Q_1
&\le&2nB_n^{-1}P(\|H\|>M_n)E(L^2)\biggl\{M_n\left(\frac\beta2B_n-M_n\right)^{-2}\\
&&\,\,\,\,\,\,\,\,\,\,\,\,\,\,\,\,\,\,\,\,+\left(\frac\beta2B_n-M_n\right)^{-1}\biggr\}\\
&\to&0
\end{eqnarray*}
as $n\longrightarrow\infty$. There is $C\in(0,\infty)$ so that
\begin{eqnarray*}
Q_2&\le&2nB_n^{-1}P(\|H\|>M_n)\left(\frac\beta2B_n-M_n\right)^{-2}E(L^2)\left\{E\|X_{n1}\|^2\right\}^{1/2}\\
&\sim&CnB_n^{-3}M_nP(\|H\|>M_n)^{3/2}\\
&\to&0\,,
\end{eqnarray*}
the equivalence in the second line following from Karamata's theorem. This shows that
$$
\lim_{n\to\infty}B_n^{-1}E\biggl\|\sum_{j=1}^nY_{nj}\one\left(\|Y_{nj}\|>\beta B_n\right)\biggr\|=0\,,
$$
and hence, showing \eqref{clt.t2.trunc} suffices for (\ref{clt.t2.eq5}).
Let
$$\sigma_{n,F}:=B_n^{-1}\sup_{f\in(B/F)^\prime,\|f\|\le1}\left[\sum_{j=1}^nf^2(T_F(u_{nj}))\right]^{1/2}\,.$$
By Theorem 4.7 in \cite{ledoux:talagrand:1991} on concentration of\\ Rademacher processes, with the median replaced by the expected value, as in page 292 of the same reference, it follows
that
\begin{eqnarray*}
&&P\left[\left|\biggl\|\sum_{j=1}^n \vep_jT_F(u_{nj})\biggr\|-E_
\vep\biggl\|\sum_{j=1}^n
\vep_jT_F(u_{nj})\biggr\|\right|>B_n\eta\right]
\le\frac{10^3}{\eta^2}E\sigma_{n,F}^2\,.
\end{eqnarray*}
Thus all
that needs to be shown is that given any $\delta>0$, there is a
choice of $\beta$ depending only on $\delta$, so that
$$\limsup_{k\to\infty}\limsup_{n\rightarrow\infty}E\sigma_{n,F_k}^2\le\delta\,.$$
Using Lemma 6.6 (page 154) in \cite{ledoux:talagrand:1991}, it
follows that for any $n,F$,
$$E\sigma_{n,F}^2\le nB_n^{-2}\sup_{f\in(B/F)^\prime,\|f\|\le1}Ef^2(T_F(u_{n1}))+8B_n^{-2}E\biggl\|\sum_{j=1}^nu_{nj}\|u_{nj}\|\biggr\|\,.$$
Clearly,
$$nB_n^{-2}\sup_{f\in(B/{F_k})^\prime,\|f\|\le1}Ef^2(T_{F_k}(u_{n1}))\le[M_n^2P(\|H\|>M_n)]^{-1}E(\|T_{{F_k}}(Y_{n1})\|^2)\,$$
which can be made as small as needed by (\ref{clt.t2.eq4}). For
the other part, note that by the contraction principle (Theorem 4.4 in \cite{ledoux:talagrand:1991}),
\begin{eqnarray*}
B_n^{-2}E\biggl\|\sum_{j=1}^nu_{nj}\|u_{nj}\|\biggr\|&\le&\beta B_n^{-1}E\biggl\|\sum_{j=1}^nu_{nj}\biggr\|\\
&\le&\beta B_n^{-1}E\biggl\|\sum_{j=1}^nY_{nj}\biggr\|\\
&=&\beta E\|Z_n\|\,.
\end{eqnarray*}
Thus,
choosing $\beta $ smaller than $\delta/(16\sup_{n\ge1}E\|Z_n\|)$
(which is positive because of (\ref{clt.t2.eq2}) ) does the trick.

 For the proof of (\ref{clt.t2.eq6}) we shall show that there is an universal constant $C>0$ so that whenever $F$ is a subspace satisfying
\begin{equation}\label{clt.t2.eq8}
\liminf_{n\to\infty}P\left[E_\vep\biggl\|\sum_{j=1}^n\vep_jT_F(Y_{nj})\biggr\|\le2B_n\delta\right]
>0\,,
\end{equation}
it follows that
\begin{equation}\label{clt.t2.eq9}
\limsup_{n\to\infty}P\left[E_\vep\biggl\|\sum_{j=1}^n\vep_jT_F(Y_{nj})\biggr\|>CB_n\delta\right]\le\delta\,.
\end{equation}
The reason that this suffices is the following. Fix $\delta>0$ and
a sequence of increasing finite-dimensional subspaces $\{F_k\}$ satisfying
\eqref{clt.closure}. Note that for all $n,k\ge1$,
\begin{eqnarray*}
&&P\left[B_n^{-1}E_ \vep\biggl\|\sum_{j=1}^n \vep_jT_{F_k}(Y_{nj})\biggr\|>2\delta\right]\\
&\le&P(\|Z_n\|>\delta)\\
&&+P\left[\left|\biggl\|\sum_{j=1}^n
\vep_jT_{F_k}(Y_{nj})\biggr\|-E_ \vep\biggl\|\sum_{j=1}^n
\vep_jT_{F_k}(Y_{nj})\biggr\|\right|>B_n\delta\right]\,.
\end{eqnarray*}
By (\ref{clt.t2.eq5}) and (\ref{clt.t2.eq1}), it follows that
$$\liminf_{k\rightarrow\infty}\liminf_{n\rightarrow\infty}P\left[E_ \vep\biggl\|\sum_{j=1}^n \vep_jT_{F_k}(Y_{nj})\biggr\|\le2B_n\delta\right]>0\,.$$
By \eqref{clt.t2.eq9}, \eqref{clt.t2.eq6} follows.

The proof of \eqref{clt.t2.eq9} uses an isoperimetric inequality;
see Theorem 1.4 (page 26) in \cite{ledoux:talagrand:1991}. Let
$$\theta:=\liminf_{n\to\infty}P\left[E_\vep\biggl\|\sum_{j=1}^n\vep_jT_F(Y_{nj})\biggr\|\le2B_n\delta\right]
>0\,.$$
In light of the isoperimetric inequality, by similar arguments as in page 291 of
\cite{ledoux:talagrand:1991}, it follows that for $k,q\ge1$,
$$\limsup_{n\to\infty}P\left[E_ \vep\biggl\|\sum_{j=1}^n \vep_jT_{F}(Y_{nj})\biggr\|>(2q+1)B_n\delta\right]$$
$$\le\left[K\left(\frac{\log(1/\theta)}k+\frac1q\right)\right]^k+P\left[B_n^{-1}\max_{j\le n}\|Y_{nj}\|>\frac\delta k\right]\,,$$
where $K$ is the universal constant in the isoperimetric
inequality. Choose $q=2K$ and $k$ to be large enough (depending
only on $\theta$) so that
$$\left[K\left(\frac{\log(1/\theta)}k+\frac1q\right)\right]^k\le\frac\delta2\,.$$
All that remains to be shown is
\begin{equation}\label{clt.t2.eq10}
\lim_{n\rightarrow\infty}P\left[B_n^{-1}\max_{j\le n}\|Y_{nj}\|>\frac\delta k\right]=0\,.
\end{equation}
Note that
$$\max_{j\le n}\|Y_{nj}\|\le\max_{j\le n}\|X_{nj}\|+\max_{j\le n}\|\tilde X_{nj}\|$$
and that
$$\max_{j\le n}\|X_{nj}\|\le M_n+\max_{j\le n}L_j\,.$$
Since $EL_1^2<\infty$, $\{n^{-1/2}\max_{j\le n}L_j\}$ is a tight family.
This shows \eqref{clt.t2.eq10} and thus establishes (\ref{clt.t2.eq9}) with $C=4q+1$ and hence completes the proof
of the direct part.

The converse is straightforward. For {\bf 1.}, note that if \eqref{clt.eq2} holds,
 by the continuous mapping theorem,
$$\lim_{n\to\infty}P(B_n^{-1}\|S_n-ES_n\|\le\epsilon)=\gamma\left(\{x\in B:\|x\|\le\epsilon\}\right)\,,$$
the right hand side being positive because in a separable Banach
space a centered Gaussian law puts positive mass on any ball with
positive radius centered at origin, see the discussion on page
60-61 in \cite{ledoux:talagrand:1991}. For proving {\bf 2.} we shall appeal to
Theorem 4.2 in \cite{acosta:gine:1979}. All that needs to be shown is
\begin{equation}\label{new}
\lim_{t\to\infty}\limsup_{n\to\infty}nE\left[\|\xi_n\|\one(\|\xi_n\|>t)\right]=0\,,
\end{equation}
where
$$
\xi_n:=U_n-E(U_n)
$$
and $U_n$ is as defined in \eqref{defU}. Note that
\begin{eqnarray*}
E\|U_n\|&\le&B_n^{-1}\left[M_n+P(\|H\|>M_n)(M_n+EL)\right]\\
&\to&0\,.
\end{eqnarray*}
Fix $t>0$. For $n$ large enough so that $\|E(U_n)\|<t/2$, it follows that
\begin{eqnarray*}
&&nE\left[\|\xi_n\|\one(\|\xi_n\|>t)\right]\\
&\le&nE\left[\|\xi_n\|\one(\|U_n\|>t/2)\right]\\
&\le&nE\left[\|U_n\|\one(\|U_n\|>t/2)\right]+nP(\|U_n\|>t/2)E\|U_n\|\,.
\end{eqnarray*}
By \eqref{new1}, the first term goes to zero. For the second term, notice that for $n$ large enough,
\begin{eqnarray*}
nP(\|U_n\|>t/2)E\|U_n\|&\le&nP(\|H\|>M_n)P\left(L>\frac t2B_n-M_n\right)E\|U_n\|\\
&=&O\left(nP(\|H\|>M_n)B_n^{-2}\right)\\
&=&o(1)\,.
\end{eqnarray*}
This shows \eqref{new} and hence completes the proof.
\end{proof}

Recall that a Banach space $B$ is said to by of type $2$ if there is $C<\infty$ so that for all $N\ge1$ and zero mean independent $B$-valued random variables $X_1,\ldots,X_N$,
$$
E\biggl\|\sum_{j=1}^n X_j\biggr\|^2\le C\sum_{j=1}^NE\|X_j\|^2\,.
$$
Banach spaces of type $2$ are nice in the sense that every random
variable $X$ taking values there with $E||X\|^2<\infty$ satisfies
the Central Limit Theorem. In fact these are the only spaces where
this is true. This is the statement of Theorem 10.5
(page 281) in \cite{ledoux:talagrand:1991}. We would like to
mention at this point that while the assumption of type $2$ is a
rather restrictive one, this is a fairly large class. For example,
every Hilbert space and $L_p$ spaces for $2\le p<\infty$ are Banach spaces of
type $2$. We show in the following result that \eqref{clt.eq2} can
be extended to these spaces.

\begin{theorem}\label{clt.t4} If $B$ is of type 2 and the model with power law tails \eqref{e:the.model} is in the hard truncation regime, then there is a Gaussian measure $\gamma$ on $B$ such that
$$B_n^{-1}(S_n-ES_n)\Rightarrow\gamma$$
The characteristic function of $\gamma$ is given by
(\ref{clt.e3}).
\end{theorem}

\begin{proof}
In view of Lemma \ref{clt.l1} and using similar arguments as in the
proof of Theorem \ref{clt.t2}, it suffices to prove that
$\{{\mathcal L}(Z_n)\}$ is relatively compact where the definition
of $Z_n$ (and $Y_{nj}$) is exactly the same as in the proof of the
latter theorem. Choose a sequence $\{F_k\}$ of finite dimensional
subspaces satisfying \eqref{clt.closure}. Since $B$ is of type 2, so is $B/F$ for any closed subspace $F$, with the type 2 constant not larger than that of $B$. Thus,
there is $C\in[0,\infty)$ so that
\begin{eqnarray*}
E\|T_{F_k}(Z_n)\|^2&\le&C[M_n^2P(\|H\|>M_n)]^{-1}E\|T_{F_k}(Y_{n1})\|^2\,.
\end{eqnarray*}
Using (\ref{clt.t2.eq4}), it follows that
$\lim_{k\rightarrow\infty}\limsup_{n\rightarrow\infty}E\|T_{F_k}(Z_n)\|^2=0$
which shows \eqref{clt.tight} and thus completes the proof.
\end{proof}

\section{Examples}\label{sec:examples} In this section, we construct a couple of examples. In Example 1, the hypotheses of Theorem \ref{clt.t2} can be verified. This helps to conclude that the result has some practical value. In Example 2, \eqref{clt.eq2} does not hold, and hence there is a need for a result like Theorem \ref{clt.t2} or Theorem \ref{clt.t4}.

\subsection*{Example 1.} Let $\{T_{jk}:j,k\ge1\}$ be i.i.d. $\bbr$-valued symmetric $\alpha$-stable (S$\alpha$S) random variables with $0<\alpha<2$, {\it i.e.}, have the following characteristic function:
$$
E[\exp i\theta T_{11}]=e^{-|\theta|^\alpha}\,.
$$
For all $j\ge1$, define the $\bbr^{\mathbb N}$-valued random variable $H_j$ as
$$
H_j:=\sum_{k=1}^\infty a_kT_{jk}e_k\,,
$$
where $(a_j)$ is a sequence of non-negative numbers satisfying
\begin{equation}\label{ex1.eq1}
\sum_{j=1}^\infty a_j^{\alpha/2}<\infty\,,
\end{equation}
and $e_k$ is the element of $\bbr^{\mathbb N}$ defined by
\begin{equation}\label{def:e}
e_k(n)=\left\{\begin{array}{ll}1,&k=n\\0,&\mbox{otherwise.}\end{array}\right.
\end{equation}
Recall that $P(|T_{11}|>x)=O(x^{-\alpha})$; see Property 1.2.15, page 16 in \cite{samorodnitsky:taqqu:1994}. This ensures that $H_1,H_2,\ldots$ are i.i.d. random variables taking values in $c_0$, the space of sequences limiting to zero, endowed with the sup norm. For that purpose, assuming that $\sum_{j=1}^\infty a_j^{\alpha}<\infty$ would have been sufficient. However, we shall need \eqref{ex1.eq1} for other reasons. It is immediate that $H_1$, and hence each $H_j$, is a $c_0$ valued symmetric $\alpha$-stable random variable. This, in particular, means that as $n\longrightarrow\infty$,
\begin{equation}\label{eq.clt}
n^{-1/\alpha}\sum_{j=1}^nH_j\Longrightarrow H_1\,.
\end{equation}
It is a well-known fact in finite dimensional spaces that the above implies
\begin{equation}\label{ex1.eq4}
P(\|H_1\|>x)\sim Cx^{-\alpha}
\end{equation}
as $x\longrightarrow\infty$ for some $C\in(0,\infty)$. However, since we could not find a reference for this on Banach spaces, we briefly sketch the argument for the sake of completeness. By Theorem 6.18, page 150 in \cite{araujo:gine:1980} it follows that there is a measure $\mu$ on $B\setminus\{0\}$ satisfying
$\mu(cD)=c^{-\alpha}\mu(D)$ for all $c>0$ and $D\subset B\setminus\{0\}$, such that,
\begin{equation}\label{eq.new2}
\lim_{n\to\infty}nP\left(n^{-1/\alpha}H_1\in A\right)=\mu(A)
\end{equation}
for all $A\subset B$ that is bounded away from the origin and $\mu(\partial A)=0$. It is also known that for all $\delta>0$, $0<\mu(\{x\in B:\|x\|\ge\delta\})<\infty$. Set
$$
A:=\{x\in B:\|x\|>1\}\,.
$$
Clearly, $\mu(\partial A)=0$. Using \eqref{eq.new2} with this $A$ implies that\\
$C:=\lim_{n\to\infty}nP(\|H_1\|>n^{1/\alpha})$ exists, and is finite and positive. Let $(x_k)$ be any sequence of positive numbers going to infinity. Set $n_k:=\lfloor x_k^\alpha\rfloor$. Observe that $x_k^\alpha P(\|H_1\|>x_k)$ is sandwiched between $n_kP(\|H_1\|>(n_k+1)^{1/\alpha})$ and $(n_k+1)P(\|H_1\|>n_k^{1/\alpha})$, and that both the bounds converge to $C$. Thus, \eqref{ex1.eq4} follows. The letter $C$ will be used to denote various such constants with possibly different definition throughout this section.

Let $(M_n)$ be a sequence such that $1\ll M_n\ll n^{1/\alpha}$. Then, the truncation of $H_j$ at level $M_n$ with $L\equiv0$ is
$$
X_{nj}:=\frac{H_j}{\|H_j\|}(\|H_j\|\wedge M_n)\,.
$$
As before, define the row sum by
$$
S_n:=\sum_{j=1}^nX_{nj}\,.
$$
We shall show that for this set up, the hypotheses of Theorem \ref{clt.t2} can be verified by purely elementary methods; the only sophisticated result that will be used is the contraction principle for finite dimensional spaces. All that needs to be shown is
\begin{equation}\label{ex1.eq2}
\liminf_{n\to\infty}P(B_n^{-1}\|S_n\|<\epsilon)>0\mbox{ for all }\epsilon>0\,,
\end{equation}
and
\begin{equation}\label{ex1.eq3}
\sup_{n\ge1}B_n^{-1}E\|S_n\|<\infty\,
\end{equation}
where
$$
B_n:=n^{1/2}M_n^{1-\alpha/2}\,.
$$
Note that in view of \eqref{ex1.eq4}, this definition of $B_n$ differs from that in the statement of Theorem \ref{clt.t2} by only a constant multiple in the limit.

The following is a sketch of how we plan to show \eqref{ex1.eq2}. Define for $K\ge0$,
\begin{eqnarray*}
S_n^{K,1}&:=&\sum_{j=1}^n\sum_{k=1}^Ka_ke_kT_{jk}\left[\one(\|H_j\|\le M_n)+\frac{M_n}{\|H_j\|}\one(\|H_j\|>M_n)\right]\,,\\
S_n^{K,2}&:=&\sum_{j=1}^n\sum_{k=K+1}^\infty a_ke_kT_{jk}\left[\one(\|H_j\|\le M_n)+\frac{M_n}{\|H_j\|}\one(\|H_j\|>M_n)\right]\,.
\end{eqnarray*}
We shall show that for all $\epsilon>0$,
\begin{equation}\label{ex1.eq6}
\sup_{K\ge1}\limsup_{n\to\infty}P\left(B_n^{-1}\|S_n^{K,1}\|>\epsilon\right)<1\,,
\end{equation}
and that
\begin{equation}\label{ex1.eq7}
\lim_{K\to\infty}\sup_{n\ge1}B_n^{-1}E\|S_n^{K,2}\|=0\,.
\end{equation}
The reason that \eqref{ex1.eq6} and \eqref{ex1.eq7} suffice for \eqref{ex1.eq2} is the following. Fix $\epsilon>0$. Note that
$$
S_n=S_n^{K,1}+S_n^{K,2}\,,
$$
and hence it follows that
$$
P(B_n^{-1}\|S_n\|>\epsilon)\le P(B_n^{-1}\|S_n^{K,1}\|>\epsilon/2)+P(B_n^{-1}\|S_n^{K,2}\|>\epsilon/2)\,.
$$
Define
$$
\delta:=1-\sup_{K\ge1}\limsup_{n\to\infty}P\left(B_n^{-1}\|S_n^{K,1}\|>\epsilon/2\right)\,,
$$
which by \eqref{ex1.eq6} is positive. Using \eqref{ex1.eq7}, choose $K$ to be large enough so that,
$$
\sup_{n\ge1}B_n^{-1}E\|S_n^{K,2}\|<\frac{\epsilon\delta}2\,.
$$
Clearly, with this choice of $K$,
$$
\limsup_{n\to\infty}P\left(B_n^{-1}\|S_n^{K,1}\|>\epsilon/2\right)\le1-\delta\,,
$$
and by the Markov inequality,
$$
\limsup_{n\to\infty}P\left(B_n^{-1}\|S_n^{K,2}\|>\epsilon/2\right)<\delta\,.
$$
This shows \eqref{ex1.eq2}.

For $n,K\ge1$, define
$$
U_{n,K}:=\sum_{j=1}^n\sum_{k=1}^Ka_ke_kT_{jk}\left[\one(a_k|T_{jk}|\le M_n)+\frac{M_n}{a_k|T_{jk}|}\one(a_k|T_{jk}|>M_n)\right]\,.
$$
We start with showing that $S_n^{K,1}$ is stochastically bounded by $U_{n,K}$, {\it i.e.}, for all $x>0$,
\begin{equation}\label{ex1.eq8}
P\left[\|S_n^{K,1}\|>x\right]\le2 P\left[\|U_{n,K}\|>x\right]\,.
\end{equation}
To that end, let $(\vep_{jk}:j,k\ge1)$ be a family of i.i.d. Rademacher random variables, independent of the family $(T_{jk}:j,k\ge1)$. Let $P_\vep$ denote the conditional probability given $(T_{jk}:j,k\ge1)$. Note that
\begin{eqnarray*}
&&P\left[\|S_n^{K,1}\|>x\right]\\
&=&P\Biggl[\Bigl\|\sum_{j=1}^n\sum_{k=1}^Ka_ke_k\vep_{jk}|T_{jk}|\left\{\one(\|H_j\|\le M_n)+\frac{M_n}{\|H_j\|}\one(\|H_j\|>M_n)\right\}\Bigr\|\\
&&\,\,\,\,>x\Biggr]\\
&=&EP_\vep\Biggl[\Bigl\|\sum_{j=1}^n\sum_{k=1}^Ka_ke_k\vep_{jk}|T_{jk}|\left\{\one(\|H_j\|\le M_n)+\frac{M_n}{\|H_j\|}\one(\|H_j\|>M_n)\right\}\Bigr\|\\
&&\,\,\,\,>x\Biggr]\,.
\end{eqnarray*}
Since the function $x\mapsto \one(x\le M)+(M/x)\one(x>M)$ is monotone non-increasing for $x\ge0$, it follows that
\begin{eqnarray}
&&|T_{jk}|\left\{\one(\|H_j\|\le M_n)+\frac{M_n}{\|H_j\|}\one(\|H_j\|>M_n)\right\}\nonumber\\
&\le&|T_{jk}|\left\{\one(a_k|T_{jk}|\le M_n)+\frac{M_n}{a_k|T_{jk}|}\one(a_k|T_{jk}|>M_n)\right\}\,.\label{ex1.eq5}
\end{eqnarray}
Using Theorem 4.4, page 95 in \cite{ledoux:talagrand:1991}, it follows that
\begin{eqnarray*}
&&P_\vep\Biggl[\Bigl\|\sum_{j=1}^n\sum_{k=1}^Ka_ke_k\vep_{jk}|T_{jk}|\left\{\one(\|H_j\|\le M_n)+\frac{M_n}{\|H_j\|}\one(\|H_j\|>M_n)\right\}\Bigr\|\\
&&\,\,\,\,>x\Biggr]\\
&\le&2P_\vep\Biggl[\Bigl\|\sum_{j=1}^n\sum_{k=1}^Ka_ke_k\vep_{jk}|T_{jk}|\Bigl\{\one(a_k|T_{jk}|\le M_n)\\
&&\,\,\,\,+\frac{M_n}{a_k|T_{jk}|}\one(a_k|T_{jk}|>M_n)\Bigr\}\Bigr\|>x\Biggr]\,.
\end{eqnarray*}
This shows \eqref{ex1.eq8}. By the result in one dimension (Theorem \ref{clt.t4} for example), it follows that for all $k\ge1$,
$$
B_n^{-1}\sum_{j=1}^na_kT_{jk}\left[\one(a_k|T_{jk}|\le M_n)+\frac{M_n}{a_k|T_{jk}|}\one(a_k|T_{jk}|>M_n)\right]\Longrightarrow N(0,a_k^{\alpha}\sigma^2)
$$
as $n\longrightarrow\infty$, where $\sigma>0$ is independent of $k$. Thus, it follows that
\begin{eqnarray*}
\lim_{n\to\infty}P\left(B_n^{-1}\|U_{n,K}\|>\epsilon/2\right)&=&1-\Pi_{k=1}^KP(|G|\le a_k^{-\alpha/2}\epsilon/2)\\
&\le&1-\Pi_{k=1}^\infty P(|G|\le a_k^{-\alpha/2}\epsilon/2)\,,
\end{eqnarray*}
where $G$ is a normal random variable with mean zero and variance $\sigma^2$. Thus, \eqref{ex1.eq6} will follows if the following is shown: for all $\eta>0$,
$$
\Pi_{k=1}^\infty P(|G|\le a_k^{-\alpha/2}\eta)>0\,.
$$
Clearly, it suffices to show that
$$
\sum_{k=1}^\infty P(|G|>a_k^{-\alpha/2}\eta)<\infty\,,
$$
which immediately follows from the Markov inequality along with \eqref{ex1.eq1}. Thus, \eqref{ex1.eq6} follows.

For showing \eqref{ex1.eq7}, note that for all $n\ge1$ and $K\ge0$,
\begin{eqnarray}
&&B_n^{-1}E\|S_n^{K,2}\|\nonumber\\
&\le&\sum_{k=K+1}^\infty E\Biggl|B_n^{-1}\sum_{j=1}^na_kT_{jk}\Bigl\{\one(\|H_j\|\le M_n)+\frac{M_n}{\|H_j\|}\one(\|H_j\|>M_n)\Bigr\}\Biggr|\nonumber\\
&\le&\sum_{k=K+1}^\infty E^{1/2}\Biggl[B_n^{-1}\sum_{j=1}^na_kT_{jk}\Bigl\{\one(\|H_j\|\le M_n)+\frac{M_n}{\|H_j\|}\one(\|H_j\|>M_n)\Bigr\}\Biggr]^2 \nonumber\\
&=&\sum_{k=K+1}^\infty B_n^{-1}n^{1/2}E^{1/2}\Biggl[a_k|T_{1k}|\Bigl\{\one(\|H_1\|\le M_n)\nonumber\\
&&\,\,\,\,+\frac{M_n}{\|H_1\|}\one(\|H_1\|>M_n)\Bigr\}\Biggr]^2 \nonumber
\end{eqnarray}
\begin{eqnarray}
&\le&\sum_{k=K+1}^\infty B_n^{-1}n^{1/2}E^{1/2}\Biggl[a_k|T_{1k}|\Bigl\{\one(a_k|T_{1k}|\le M_n)\label{ex1.eq9}\\
&&\,\,\,\,+\frac{M_n}{a_k|T_{1k}|}\one(a_k|T_{1k}|>M_n)\Bigr\}\Biggr]^2 \nonumber\\
&\le&C\sum_{k=K+1}^\infty a_k^{\alpha/2}\label{ex1.eq10}
\end{eqnarray}
for some $C<\infty$ independent of $n$ and $K$, where \eqref{ex1.eq5} has been used for \eqref{ex1.eq9}, and \eqref{ex1.eq10} follows by Karamata Theorem, the estimation being similar to that leading to \eqref{clt.new}. This, in view of \eqref{ex1.eq1}, shows \eqref{ex1.eq7}. Thus, \eqref{ex1.eq2} follows. Also, using \eqref{ex1.eq10} for $K=0$, \eqref{ex1.eq3} follows. Thus, the hypotheses of Theorem \ref{clt.t2} are satisfied.

\subsection*{Example 2}Fix $1<p<2$. We first construct a bounded symmetric random variable $X$ taking values in $c_0$ (the space of sequences limiting to zero, equipped with the sup norm) so that $n^{-1/p}\sum_{i=1}^nX_i$ does not converge to zero in probability, where $X_1,X_2,\ldots$ are i.i.d. copies of $X$. Let $(\vep_j:j\ge1)$ be a sequence of i.i.d. Rademacher random variables. We shall use the fact that there exists $K\in(0,\infty)$ so that
\begin{equation}\label{ex2.eq1}
P\left(\sum_{i=1}^n\vep_i>t\right)\ge\exp\left(-Kt^2/n\right)
\end{equation}
for all $n\ge1$ and $t>0$ such that $n^{1/2}K\le t\le K^{-1}n$. This follows from (4.2) on page 90 in \cite{ledoux:talagrand:1991}. Define
$$
X:=\sum_{j=1}^\infty a_j\vep_je_j\,,
$$
where
$$
a_j:=K\{\log(j\vee2)\}^{(1-p)/2},\,j\ge1\,,
$$
$K$ is the constant in \eqref{ex2.eq1} and $e_j$ is as defined in \eqref{def:e}. Clearly, $X$ is a $c_0$ valued symmetric bounded random variable. Let $X_1,X_2,\ldots$ denote i.i.d. copies of $X$.
Note that for $n\ge1$,
\begin{eqnarray*}
P\left(n^{-1/p}\Bigl\|\sum_{k=1}^nX_k\Bigr\|>1\right)&=&1-\Pi_{j=1}^\infty P\left(\Bigl|\sum_{k=1}^n\vep_k\Bigr|\le n^{1/p}a_j^{-1}\right)\,.
\end{eqnarray*}
Thus, for proving that $n^{-1/p}\sum_{k=1}^nX_k$ does not converge to zero in probability, it suffices to show that
\begin{equation}\label{ex2.eq2}
\limsup_{n\to\infty}\Pi_{j=1}^\infty P\left(\Bigl|\sum_{k=1}^n\vep_k\Bigr|\le n^{1/p}a_j^{-1}\right)<1\,.
\end{equation}
To that aim, define
$$
l_n:=\left\lfloor\exp\left(n^{2/p}\right)\right\rfloor,\,n\ge1\,,
$$
and note that
\begin{eqnarray*}
\Pi_{j=1}^\infty P\left(\Bigl|\sum_{k=1}^n\vep_k\Bigr|\le n^{1/p}a_j^{-1}\right)
&\le&\Pi_{j=1}^{l_n} P\left(\Bigl|\sum_{k=1}^n\vep_k\Bigr|\le n^{1/p}a_j^{-1}\right)\\
&\le&P\left(\Bigl|\sum_{k=1}^n\vep_k\Bigr|\le n^{1/p}a_{l_n}^{-1}\right)^{l_n}\,.
\end{eqnarray*}
Note that
\begin{eqnarray*}
n^{1/p}a_{l_n}^{-1}&=&K^{-1}n^{1/p}(\log l_n)^{-(1-p)/2}\\
&\le&K^{-1}n\,.
\end{eqnarray*}
Also, it is easy to see that as $n\longrightarrow\infty$,
\begin{equation}\label{ex2.eq3}
\log l_n\sim n^{2/p}\,.
\end{equation}
Thus, $n^{1/p}a_{l_n}^{-1}\gg n^{1/2}$. For $n$ large enough, an appeal to \eqref{ex2.eq1} shows that
\begin{eqnarray*}
P\left(\Bigl|\sum_{k=1}^n\vep_k\Bigr|\le n^{1/p}a_{l_n}^{-1}\right)
&\le&1-\exp\left(-Kn^{2/p-1}a_{l_n}^{-2}\right)\,.
\end{eqnarray*}
Using \eqref{ex2.eq3}, it follows that
$$
n^{2/p-1}a_{l_n}^{-2}=o(\log l_n)\,.
$$
Thus, for $n$ large enough it holds that
$$
Kn^{2/p-1}a_{l_n}^{-2}\le\log l_n\,,
$$
and hence for such a $n$,
$$
\exp\left(-Kn^{2/p-1}a_{l_n}^{-2}\right)\ge\frac1{l_n}\,.
$$
What we have shown can be summed up as that for $n$ large enough,
$$
\Pi_{j=1}^\infty P\left(\Bigl|\sum_{k=1}^n\vep_k\Bigr|\le n^{1/p}a_j^{-1}\right)\le\left(1-\frac1{l_n}\right)^{l_n}\,.
$$
Thus, \eqref{ex2.eq2} follows.

Fix $x\in B\setminus\{0\}$ and define
\begin{equation}\label{defY}
Y:=X\one(U=0)+xS\one(U=1)
\end{equation}
where $S$ is a $\bbr$-valued (symmetric) Cauchy random variable and $U$ is a Bernoulli$(1/2)$ random variable such that $X,S,U$ are all independent.
We start with showing that $Y$ is in the domain of attraction of an
$1$-stable law on $B$.
Let $((X_i,S_i,U_i):i\ge1)$ denote i.i.d. copies of $(X,S,U)$. Since $X$ has zero mean, by Theorem 9.21 in \cite{ledoux:talagrand:1991}, it follows that
$$n^{-1}\sum_{i=1}^nX_i\prob0\,.$$
We shall show by an application of the contraction principle (Theorem 4.4 in \cite{ledoux:talagrand:1991}) that
\begin{equation}\label{ex2.eq4}
n^{-1}\sum_{i=1}^nX_i\one(U_i=0)\prob0\,.
\end{equation}
Let $(\vep_j:j\ge1)$ be a sequence of i.i.d. Rademacher random variables independent of $((X_i,S_i,U_i):i\ge1)$. Let $P_\vep$ denote the conditional probability given $((X_i,S_i,U_i):i\ge1)$. Thus for all $u>0$,
\begin{eqnarray*}
P\left(\Bigl\|\sum_{i=1}^nX_i\one(U_i=0)\Bigr\|>u\right)
&=&EP_\vep\left(\Bigl\|\sum_{i=1}^n\vep_iX_i\one(U_i=0)\Bigr\|>u\right)\\
&\le&2EP_\vep\left(\Bigl\|\sum_{i=1}^n\vep_iX_i\Bigr\|>u\right)\\
&=&2P\left(\Bigl\|\sum_{i=1}^nX_i\Bigr\|>u\right)\,,
\end{eqnarray*}
the inequality following by the contraction principle. This shows \eqref{ex2.eq4}.
By Theorem 3 on page 580 in \cite{feller:1971}, it follows that
$$
n^{-1}\sum_{i=1}^nS_i\one(U_i=1)\Longrightarrow Z
$$
for some Cauchy random variable $Z$. Thus, it is immediate that
\begin{equation}\label{eq.new3}
n^{-1}\sum_{i=1}^nY_i\Longrightarrow xZ\,,
\end{equation}
where $Y_1,Y_2,\ldots$ denote i.i.d.
copies of $Y$.

For a positive number $M_n$,
$$Y_{ni}:=Y_i\one\left(\|Y_i\|\le M_n\right)+M_n\frac{Y_i}{\|Y_i\|}\one\left(\|Y_i\|>M_n\right)$$
is the truncation of $Y_i$ to the ball of radius $M_n$, as defined in \eqref{e:the.model} with $L$ identically equal to zero. Let
$$S_n:=\sum_{i=1}^nY_{ni}\,.$$
We will show $n^{-1/p}S_n$ does not converge to $0$ in probability whenever
$M_n\longrightarrow\infty$.
By arguments similar to those leading to \eqref{ex2.eq4}, it follows that for $u>0$,
\begin{eqnarray*}
P\left(\Bigl\|\sum_{i=1}^nY_i\one\left(\|Y_i\|\le M_n\right)\Bigr\|>u\right)
&\le&2P(\|S_n\|>u)\,.
\end{eqnarray*}
Note that since $X$ is bounded and $M_n$ goes to infinity, for $n$ large enough,
$$
Y\one(\|Y\|\le M_n)=X\one(U=0)+xS\one(U=1)\one(|S|\le M_n/\|x\|)\,.
$$
Observing that if $(\vep_1,\vep_2)$ are i.i.d. Rademacher random variables independent of $(X,S,U)$, then
$$
X\one(U=0)+xS\one(U=1)\one(|S|\le M_n/\|x\|)
$$
$$
\eid\vep_1X\one(U=0)+x\vep_2|S|\one(U=1)\one(|S|\le M_n/\|x\|)\,,
$$
exactly same arguments as before will show that for $n$ large enough and $u>0$,
$$
P\left(\Bigl\|\sum_{i=1}^nX_i\one(U_i=0)\Bigr\|>u\right)\le2P\left(\Bigl\|\sum_{i=1}^nY_i\one\left(\|Y_i\|\le M_n\right)\Bigr\|>u\right)\,.
$$
The above can be summarized as that there exists $N<\infty$ so that
$$
P\left(\Bigl\|\sum_{i=1}^nX_i\one(U_i=0)\Bigr\|>u\right)\le4(P(\|S_n\|>u)
$$
for all $n\ge N$ and $u>0$.

Denote
$$N_n:=n-\sum_{i=1}^nU_i\,,$$
and the conditional probability given $U_1,U_2,\ldots$ by $P_U$. Note that
\begin{eqnarray*}
&&P\left(\Bigl\|\sum_{i=1}^nX_i\one(U_i=0)\Bigr\|>u\right)\\
&\ge&P\left(\Bigl\|\sum_{i=1}^nX_i\one(U_i=0)\Bigr\|>u,N_n>n/3\right)\\
&=&E\left[P_U\left(\Bigl\|\sum_{i=1}^{N_n}X_i\Bigr\|>u\right)\one(N_n>n/3)\right]\,.
\end{eqnarray*}
Another application of the contraction principle shows that on the set $\{N_n>n/3\}$,
$$
P_U\left(\Bigl\|\sum_{i=1}^{N_n}X_i\Bigr\|>u\right)\ge\frac12P_U\left(\Bigl\|\sum_{i=1}^{\lceil n/3\rceil }X_i\Bigr\|>u\right)=\frac12P\left(\Bigl\|\sum_{i=1}^{\lceil n/3\rceil }X_i\Bigr\|>u\right)\,.
$$
Thus, it follows that
$$
P\left(\Bigl\|\sum_{i=1}^nX_i\one(U_i=0)\Bigr\|>u\right)\ge\frac12P\left(\Bigl\|\sum_{i=1}^{\lceil n/3\rceil }X_i\Bigr\|>u\right)P(N_n>n/3)\,.
$$
All the above calculations put together shows
$$
P\left(\Bigl\|\sum_{i=1}^{\lceil n/3\rceil }X_i\Bigr\|>u\right)=O\left(P(\|S_n\|>u)\right)
$$
uniformly in $u$.
Since $n^{-1/p}\sum_{i=1}^nX_i$ does not converge to zero in probability, it follows that $n^{-1/p}S_n$ does not converge to $0$ in probability either.

The above calculations can be used to construct an example where \eqref{clt.eq2} does not hold, in the following way. Fix $1<p<2$ and a sequence $(M_n)$ satisfying $1\ll M_n\ll n^{2/p-1}$. Define $Y$ by \eqref{defY}. The argument that leads to \eqref{ex1.eq4} from \eqref{eq.clt}, applied to $Y$ helps us conclude from \eqref{eq.new3} that
$
P(\|Y\|>x)\sim Cx^{-1}
$
as $x\longrightarrow\infty$, for some $C\in(0,\infty)$. Since $2/p-1<1$, it follows that $M_n\ll n$, which is a  restatement of
\begin{eqnarray*}
\lim_{n\to\infty}nP(\|Y\|>M_n)=\infty\,.
\end{eqnarray*}
Thus, the assumption of hard truncation is satisfied. Set $L\equiv0$, $Y_{ni}$ to be the truncation of $Y_i$ at level $M_n$, $S_n$ to be the row sum of the triangular array $\{Y_{ni}:1\le i\le n\}$ and
$$
B_n:=\left[nM_n^2P(\|Y\|>M_n)\right]^{1/2}\,.
$$
Thus,
$$
B_n^2=O(nM_n)=o\left(n^{2/p}\right)\,.
$$
This shows that $B_n^{-1}S_n$ does not converge weakly, for otherwise, $n^{-1/p}S_n$ would converge to zero in probability. Thus, \eqref{clt.eq2} does not hold.

This is an example where the claim of Theorem \ref{clt.t4} does not hold. The space $c_0$ is not of Rademacher type $p$ for all $p>1$. Hence it was possible to construct a zero mean random variable with finite $p$-th moment, that does not satisfy the law of large numbers with rate $n^{1/p}$.

\section{Acknowledgment} The author is immensely grateful to his adviser Gennady Samorodnitsky for many helpful discussions, and to Parthanil Roy for his suggestions on the layout. The author also thanks an anonymous referee for some suggestions that significantly improved the paper.

\end{document}